\documentclass[preprint,a4paper]{article}

\usepackage{amssymb}

\newtheorem{problem}{\bf Problem}

\title{Number systems and combinatorial problems}
\author{Krasimir Yordzhev}
\date{\empty}

\begin{document}
\maketitle

\begin{center}
South-West University, Faculty of Mathematics and Natural Sciences\\
66 Ivan Mihailov Str, 2700 Blagoevgrad, Bulgaria \\
E-mail: yordzhev@swu.bg
\end{center}

\begin{abstract}
 The present work has been designed for students in secondary school and their teachers in mathematics. We will show how with the help of our knowledge of number systems we can solve problems from other fields of mathematics – for example in combinatorial analysis and most of all when proving some combinatorial identities. To demonstrate discussed in this article method we have chosen several suitable mathematical tasks.

\textbf{Keywords}: number system; binary notation;  $n$-digit  $p$-ary number; combinatorial identity

\textbf{2010 MSC}:   05A19; 97K20
\end{abstract}

With the development of computer technology we can observe general application of different mathematical models in all spheres of science, including humanitarian disciplines. The leader in this regard is computer science. It is impossible for a person to be engaged in programming and information technologies without any solid mathematical knowledge and skills. Here we are going to discuss the other side of the problem – how a student interested in computers and information technologies can consider some theoretical mathematical statements.

The great application of different number systems and most of all the binary notation in computing machinery and in programming is well-known. Here we are not going to consider these problems, but we will try to show how with the help of our knowledge of number systems we can solve problems from other fields of mathematics – for example in combinatorics and most of all when proving some combinatorial identities.

Let us first remind the reader of what the symbol   $n\choose k$, which plays an important role in combinatorial analysis and is related to the development of binomial formula, means. This is the number of all combinations of $n$  elements, class  $k$, i.e. $n\choose k$  shows in how many ways we can choose   $k$ objects belonging to a  set with  $n$ element. Probably most of the students know how to express $n\choose k$   by $n$  and  $ k$, but let us forget about the formula of this representation for a moment. In the solutions to the suggested problems we can skip this formula.

Under  $n$-digit  $p$-ary number we will perceive a number written in notation with the base $p$ in exactly $n$  positions with the help of the digits  $\{0,1,...,p-1\}$. We hope that the reader is familiar with the basic concepts related to notations and that he/she can solve at least elementary problems from this field. If a number can be written with less than $n$  digits in a  $p$-ary notation then the  $n$-digit  $p$-ary number is completed from the left with the necessary number of zeros. In particular when $p=2$  we obtain the binary notation, where each  $n$-digit binary number is written with exactly $n$  digits 0 or 1. Then $n\choose k$  can be interpreted as the number of all  $n$-digit binary numbers with exactly   $k$ 1's and $(n-k)$ 0's. But this number  equals the number of all  $n$-digit binary numbers with exactly  $k$ 0's  and $(n-k)$ 1's. In this way we have reached the well-known equation

$${n\choose k} = {n\choose n-k} .$$

The next statement is a well-known fact, for which we hope that the reader can point out more than one proof. With this problem we are willing to illustrate the methods of proof, which are the object of consideration in the present article.

\begin{problem}\label{z1}
Prove the following identity:
$$\sum_{k=0}^{n}{n\choose k}={n\choose 0} + {n\choose 1} +\cdots+ {n\choose n} = 2^n $$
\end{problem}

Proof. Let   $S=\displaystyle \sum_{k=0}^n {n\choose k}$. Obviously $S$  equals the number of all subsets (including the empty one) of the set  $X=\{x_1,x_2, ... ,x_n\}$. Each subset  $A\subseteq X$ can be coded with the  $n$-digit binary number  $a=\alpha_1 \alpha_2 \ldots \alpha_n$, where  $\alpha_i=1$, if $x_i\in A$  and  $\alpha_i=0$, if  $x_i\notin A$, $i=1,2,\ldots,n$. Conversely for every binary number  $a=\alpha_1 \alpha_2 \ldots \alpha_n$  corresponds to a subset $A\subseteq X$  in the following way: $x_i\in A$  if and only if $\alpha_i=1$.

Therefore there is one-to-one correspondence between the all subsets of $X$  and all  $n$-digit binary numbers     $\displaystyle \underbrace{00...0}_n,$ $\underbrace{00...1}_n,$ $\ldots ,$ $\underbrace{11...1}_n $, whose decimal notation is respectively  $0,1,2,...,(2^n-1)$, and these numbers are exactly $2^n$  in number.

\bigskip
The next problem is a generalization of problem 1. We are leaving its solution to the reader for his/her independent work.

\begin{problem}\label{z2}
 Let us have the set  $X=\lbrace x_1,x_2,...,x_n \rbrace $, composed of  $n$ objects. Each object can be coloured in a random way, while the set of colours is  $\lbrace c_0,c_1,...,c_{p-1} \rbrace $, i.e. the colours are $p$  in number. With ${n\brack k}_{p,i}$  we are denoting the number of all colourings of the elements of the set  $X$, so that exactly $k$  in number objects are painted in the colour  $c_i$. (Obviously  ${n\brack k}_{p,0}={n\brack k}_{p,1}=\cdots={n\brack k}_{p,(p-1)}$.) Then for each   $i=0,1,...,(p-1)$, prove the following identity

$${n\brack 0}_{p,i}+{n\brack 1}_{p,i}+\cdots +{n\brack n}_{p,i}=p^n$$
\end{problem}

The proof to problem 2 is analogical to the proof to problem 1, only instead of a binary notation we will have to use the notation with the base $p$.

\bigskip

The binary functions with  $n$ variables, where   $n$ is a positive integer, are greatly applied in computer science. These are functions from the type  $f(x_1 ,x_2 ,\ldots ,x_n )$, whose variables  $x_i$, $i=1,2, \ldots ,n$, as well as the function $f$ take the value only of the set  $\{ 0,1\}$. Likewise a   $k$-valued function with $n$  variables is called a function from the type  $f(x_1 ,x_2 ,\ldots ,x_n )$, whose variables $x_i$, $i=1,2, \ldots ,n$, as well as the function $f$ take value only of  the set  $\{ 0,1, \ldots , k-1\}$, where  $k\ge 2$ is a positive integer. Our knowledge in number systems will help us solve the next two problems, which are basic tasks in Boolean algebra and   $k$-valued logic.

\begin{problem}
Find the number of all binary functions with $n$  variables, where $n$ is a positive integer.
\end{problem}

Solution. Each binary function with $n$  variables can be presented with the help of a table consisting of  $n+1$ columns.
$$
\begin{array}{cccc|c}
  x_1 & x_2 & \cdots & x_n & f(x_1 ,x_1 ,\ldots ,x_n ) \\
  \hline\hline
  0 & 0 & \cdots & 0 & f(0,0,\ldots ,0)=\beta_1 \\
  0 & 0 & \cdots & 1 & f(0,0,\ldots ,1)=\beta_2 \\
  \vdots & \vdots & \cdots & \vdots & \vdots \\
  \alpha_{i_1} & \alpha_{i_2} & \cdots & \alpha_{i_n} &  f(\alpha_{i_1} ,\alpha_{i_2} , \ldots , \alpha_{i_n}) =\beta_i \\
  \vdots & \vdots & \cdots & \vdots & \vdots \\
  1 & 1 & \cdots & 1 &  f(1,1,\ldots ,1)=\beta_{2^n} \\
\end{array}
$$
In the first $n$  columns of the table we put all ordered  $n$-tuples from the type  $\alpha_{i_1} ,\alpha_{i_2} ,\ldots ,\alpha_{i_n}$, where  $\alpha_{i_j} \in \{ 0,1\}$,  $j=1,2,\ldots n$. In the different rows we put different  $n$-tuples and these are  $\underbrace{0,0,\ldots ,0}_n$; $\underbrace{0,0,\ldots ,1}_n$; $\ldots$ ; $\underbrace{1,1,\ldots ,1}_n$. Therefore in the first  $n$ columns of the table all integers belonging to the interval from  $0=\underbrace{00\ldots 0}_n$ to  $2^n -1=\underbrace{11\ldots 1}_n$ are written in binary notation.

Therefore the number of the rows in the table equals   $2^n$.

In the last column of the table we put the corresponding values of the function, i.e. if in the   $i$-th row of the chart in the first $n$  columns we have put the  $n$-tuple  $\alpha_{i_1} ,\alpha_{i_2} ,\ldots ,\alpha_{i_n}$, then in the   $(n+1)$-th column in the  $i$-th row the value of the function $f(\alpha_{i_1} ,\alpha_{i_2} ,\ldots ,\alpha_{i_n} )=\beta_i \in \{ 0,1\}$ is written. Two binary functions differ from each other if and only if there is at least one integer  $i$, such that $1\le i\le 2^n$  and the values  $\beta_i$  written in  $i$  row $n+1$  column of the corresponding table to each of the two functions differ from each other. Therefore the number of all binary functions with $n$  variables is equal to the number of all  $2^n$-tuples  $\beta_1 ,\beta_2 ,\ldots ,\beta_{2^n}$, $\beta_i \in\{ 0, 1\}$, $i=1,2,\ldots ,2^n$. Therefore the number of all binary functions is equal to the number of all integers belonging to the interval from $0=\underbrace{00\ldots 0}_{2^n} $  to  $2^n -1 =\underbrace{11\ldots 1}_{2^n} $  (written in a binary notation), i.e. the number of the binary functions with   $n$ variables equals  $2^{2^n}$.

\bigskip

Likewise we can solve the following problem:

\begin{problem}
Find the number of all  $k$-valued functions with $n$  variables, where $k$ and $n$ are positive integers, $k\ge 2$.
\end{problem}

By using the properties of the number system with a given base, it is not difficult to solve the following problem.

\begin{problem}
Consider the equation $xy=a^p$, where   $a$   and  $p$ are natural numbers and $a$  can be presented as a product of   $n$ different primes. Find the number of the solutions to the above equation in natural numbers.
\end{problem}

 Solution. Let  $a=a_1a_2\ldots a_n$, where $a_i$  is a prime number,  $i=1,2,\ldots,n$. Obviously  all solutions to the given equation are presented in the form of  $x=a_1^{\beta_1}a_2^{\beta_2}\ldots a_n^{\beta_n}$, $y=a_1^{p-\beta_1}a_2^{p-\beta_2}\ldots a_n^{p-\beta_n}$, where  $\beta_i\in\lbrace 0,1,...,p\rbrace$, $i=1,2,...,n$. Therefore the  $n$-digit  $(p+1)$ -ary number $\beta_1\beta_2\ldots\beta_n$  corresponds to every solution to the given equation. Conversely to every  $n$-digit  $(p+1)$-ary number $\beta_1\beta_2\ldots\beta_n$ corresponds exactly one solution. Bearing in mind that the number of all  $n$-digit  $(p+1)$-ary numbers is equal to  $(p+1)^n$, it follows that the number of the solutions to the given equation in natural numbers is equal to  $(p+1)^n$.

\bigskip

To the fans of mathematical puzzles we will propose the following

\begin{problem}
In how many different ways, moving from letter to letter, can we read the word “mathematics” on the figure shown below? On the figure the letters in bold show a possible ''walk''.
$$
\begin{array}{ccccccccccc}
  \textbf{m} & \textbf{a} & \textbf{t} & h & e & m & a & t & i & c & s \\
  a & t & \textbf{h} & \textbf{e} & m & a & t & i & c & s &    \\
  t & h & e & \textbf{m} & \textbf{a} & \textbf{t} & i & c & s &   &   \\
  h & e & m & a & t & \textbf{i} & c & s &   &   &   \\
  e & m & a & t & i & \textbf{c} & \textbf{s} &   &   &   &   \\
  m & a & t & i & c & s &   &   &   &   &   \\
  a & t & i & c & s &   &   &   &   &   &   \\
  t & i & c & s &   &   &   &   &   &   &   \\
  i & c & s &   &   &   &   &   &   &   &   \\
  c & s &   &   &   &   &   &   &   &   &   \\
  s &   &   &   &   &   &   &   &   &   &
\end{array}
$$
\end{problem}

Solution. The condition of the problem clarifies that standing on a letter, unless it is the last one, then the next letter from the ''walk'' is either on its right side, or below it. Since the word “mathematics” consists of eleven letters, and bearing in mind that after the last letter ''s'' there is nowhere to go, it is not difficult to notice that every ''walk'' can be coded with ten-digit binary number $\alpha_1\alpha_2\ldots\alpha_{10}$,  $\alpha_i\in\{0,1\}$, $i=1,2,\ldots ,10$  in the following way:  $\alpha_i=1$, if from the   $i$-th letter of the word mathematics we have started to its right and  $\alpha_i=0$, if we have started downwards. Conversely, apparently to every ten-digit binary number corresponds a ''walk'' according to the described rule. Therefore, the number of all 'walks'' is equal to  $2^{10}=1024$.

\bigskip
And now we are going to consider a problem, which was given at the 28th International Mathematical Olympiad. We hope that after everything we have said so far you will find it too easy. First we are going to recall the following

\textbf{Definition}: A \emph{permutation}  $f$ of the set  $[n]=\{1,2,...,n\}$  is  every ordered   $n$-tuple  $a_1a_2...a_n$, where $a_i\in [n]$  and $a_i \ne a_j$ when  $i \ne j$.  It is well-known that the number of all permutations from $n$  elements is equal to  $n!=1\cdot 2 \cdot 3 \cdots (n-1)n$. The integer  $r\in [n]$ is called fixed point of the permutation   $f$, if  $a_r=r$.

\begin{problem}
Let  $p_n(k)$  be the number of the permutations of the set   $[n]=\{ 1,2,\ldots ,n\}$, which have exactly  $k$ fixed points. Prove the equality:
$$\sum_{k=0}^n kp_n(k)=n!$$
\end{problem}

Proof. To every permutation of $[n]$  we assign the $n$-digit binary number  $\alpha_1\alpha_2\ldots\alpha_n$, so that   $\alpha_i=1$, if $i$  is a fixed point for the considered permutation and $\alpha_i=0$  otherwise ($i=1,2,\ldots ,n$). Let $A$  be the multiset of those  $n$-digit binary numbers, to which a permutation according to the above rule is assigned. (We must remark that in $A$  there are repetitive identical elements, corresponding to different permutations, while contrariwise not all  $n$-digit binary numbers belong to  $A$. For example, there are not any permutations with exactly $(n-1)$  fixed points.) Apparently $p_n(k)$  can be viewed as the number of those elements of   $A$, which have exactly  $k$ 1's. Let us denote the total number of the units in $A$  with  $q$. Then obviously  $\displaystyle q=\sum_{k=0}^n kp_n(k)$. Contrariwise, bearing in mind that the number of the permutations of $(n-1)$  objects is equal to  $(n-1)!$, it is not difficult to see that for each $i\in \{1,2,...,n\}$  there are exactly $(n-1)!$  permutations, for which the point $i$  is fixed, i.e. for which  $\alpha_i=1$. Counting the units in  $A$ in this way, we come to the conclusion that  $q=n (n-1)!=n!$. Therefore   $\displaystyle \sum_{k=0}^n kp_n(k)=n!$.

\bigskip
Finally we will give you a chance to test your strength with proving some combinatorial identities. We advise you to try to use the method exposed in the article.

\begin{problem}
If  $m$, $n$, $p$ and  $k$ are natural numbers, prove the identities:

\begin{description}
  \item[a)]  $\displaystyle{n\choose p} = {{n-1}\choose p} + {{n-1}\choose {p-1}}$

 \item[b)] $\displaystyle {{n-p}\choose {m-p}} = {{n-p}\choose {n-m}}$

  \item[c)] $\displaystyle {{n+p}\choose {n-m}} = {{n+p}\choose {m+p}}$

\item[d)] $\displaystyle {{n+m}\choose p} = \sum_{k=0}^p {n\choose k} {m\choose{p-k}}$

\item[e)] $\displaystyle {{2n}\choose n} =\sum_{k=0}^n {n\choose k}^2$

\item[f)] $\displaystyle {n\choose m} {{n-m}\choose p} = {n\choose p} {{n-p}\choose m} = {n\choose {m+p}} {{m+p}\choose m}$

\item[g)] $\displaystyle {n\choose m} {m\choose p} = {n\choose p} {{n-p}\choose {m-p}} = {n\choose {m-p}} {{n-m+p}\choose p}$

\item[h)] For each $i=0,1,\ldots ,p-1$ is satisfied  $\displaystyle {n\brack k}_{p,i} = {n\choose k}(p-1)^{n-k}$, where $\displaystyle {n\brack k}_{p,i}$  is the symbol defined in problem \ref{z2}.

\item[i)] $\displaystyle \sum_{k=0}^n{{n\choose k}(p-1)^{n-k}} = p^k$

\item[j)] $\displaystyle \sum_{k=m}^n {n\choose k} {k\choose m} = 2^{n-m} {n\choose m}$

\item[k)] $\displaystyle \sum_{k=0}^n {{2n+1}\choose k} = 2^{2n}$

\item[l)] $\displaystyle \sum_{k=0}^n k {{n}\choose k} = n2^{n-1}$;

\item[m)] $\displaystyle \sum_{p=0}^n k {{k+p}\choose p} =
{{k+n+1}\choose n} $;

\item[n)] $\displaystyle \sum_{k=0}^n  {{2k}\choose k}{{2n-2k}\choose
{n-k}} =4^n$.
\end{description}
\end{problem}

Instructions:

a) Consider the set of the all $n$-digit binary numbers as a union of two nonintersecting subsets – the subsets of the  $n$-digit binary numbers ending respectively in 0 or 1.

b) If in a   $(n-p)$-digit binary number we have exactly $(m-p)$  1's, then in this number we will have exactly $(n-p)-(m-p) = (n-m)$  0's.

c) Like b).

d) Calculate the number of all  $(n+m)$-digit binary numbers written with exactly $p$  in number 1's, while reporting the number of the 1's in the first $n$  positions.

e) Look at d) with   $n=m=p$.

f) Calculate in three different ways the number of the  $n$-digit ternary numbers, written with exactly  $m$  in number 2's and $p$  in number 1's.

g) We reason just like in f).

h) ${n\brack k}_{p,i}$  is equal to the number of all  $n$-digit  $p$-ary numbers, written with exactly $k$  digits i.

i) Follows immediately from h) and problem \ref{z2}.

j) Calculate the number of all  $n$-digit ternary numbers, written with the help of exactly $m$  2's. On the left side of the equality  $k$  it can be interpreted as the number of the digits different from 0 when writing a certain ternary number.

k) The number of all  $(2n+1)$-digit binary numbers with $k$  zeros and $2n+1-k$ units is equal to the number of all  $(2n+1)$-digit binary numbers with $2n+1-k$ zeros and $k$  units, where  $k=0,1,2,...,n$, i.e. the expression on the left side is equal to exactly half the number of all  $(2n+1)$-digit binary numbers and look at problem \ref{z1}.

Think independently on identities l), m) and n).











\end{document}